\documentclass[a4paper,reqno,11pt,final]{amsart}
\usepackage{amsthm,amsmath,amssymb}
\usepackage{mathabx} 

\usepackage[dvips]{graphicx, color}
\usepackage{graphicx}
\usepackage[title]{appendix}
 \usepackage{color}

\usepackage{color}
\usepackage{xcolor}
\usepackage{bbm}
\usepackage{scalerel}
\usepackage{stackengine}

\newtheorem{thm}  {Theorem}
\newtheorem*{theorem*}{Theorem}

\newtheorem*{theorem-non}  {Theorem}

\newtheorem{prop}[thm]{Proposition}

\newtheorem*{condition-non}  {Condition}
\newtheorem{ex}[]{Example}
\newtheorem{claim}[]{Claim}

\newtheorem{cor}[thm]{Corollary}
\newtheorem{conj}[thm]{Conjecture}
\newtheorem{lem}[thm]{Lemma}

\newtheorem{defn}[thm]{Definition}

\def\elem{\end{lem}}
\def\eth{\end{thm}}
\def\eco{\end{cor}}
\def\bco{\begin{cor}}
\def\bth{\begin{thm}}
\def\blem{\begin{lem}}
\def\epro{\end{prop}}

\def\bpro{\begin{prop}}

\newcommand{\cE}  {{\mathcal E}}

\newcommand{\cL}  {{\mathcal L}}

\newcommand{\cV}  {{\mathcal V}}

\newcommand{\E}  {{\Bbb E}}
\newcommand{\F}  {{\Bbb F}}

\renewcommand{\O}  {{\Bbb O}}
\renewcommand{\P}  {{\Bbb P}\hspace{.05em}}
\newcommand{\Q}  {{\Bbb Q}\hspace{.06em}}
\newcommand{\R}  {{\Bbb R}}

\newcommand{\Z}  {{\Bbb Z}}

\def\de{\delta}

\def\O{\mathcal O}

\def\={\:=\:}  \def\+{\,+\,}

\def\a{\alpha} \def\b{{\beta}}  \def\ba{\overline\a}    
   \def\r{\bold r}

\def\be{\begin{equation}}   \def\ee{\end{equation}}
\def\bes{\begin{equation*}}   \def\ees{\end{equation*}}
\def\ba{\begin{aligned}}   \def\ea{\end{aligned}}
\def\bc{\begin{cases}}   \def\ec{\end{cases}}
\def\bp{\begin{proof}}   \def\ep{\end{proof}}

\newcommand{\Res}  {\mathrm{Res}}

\newcommand{\Aut}  {\mathrm{Aut}}

\newcommand{\sep}  {\mathrm{\sep}}

\def\SL{\mathrm{SL}}

\def\qqan{\qquad\mathrm{and}\qquad}

\def\qfo{\quad(\forall}

\def\smm{\smallsetminus}

\def\SL{\mathrm{SL}}

\def\Lra{\Longrightarrow}

\def\De{\Delta}

\def\bbm1{\mathbbm_1}

\def\De{\Delta}

\def\wh{\widehat}
\def\be{\begin{equation}}   \def\ee{\end{equation}}
\def\bes{\begin{equation*}}   \def\ees{\end{equation*}}

\def\bea{\begin{equation}\begin{aligned}}   
\def\eea{\end{aligned}\end{equation}}

\def\r{\right}
\def\bm{\begin{matrix}}
\def\em{\end{matrix}}
\def\bpm{\begin{pmatrix}}
\def\epm{\end{pmatrix}}

\def\Pic{\mathrm{Pic}}

\def\h2m{\hskip 2.0cm}

\def\ed{\end{document}}

\def\sep{\mathrm{sep}}

\def\Res{\mathrm{Res}}

\def\Th{\Theta}
\def\th{\theta}

\usepackage{listings}
\usepackage{float}
\usepackage{color}
\usepackage{bbm}

\usepackage{tikz}
\usetikzlibrary{quantikz}
\def\qw{\qquad{\rm where}\quad}

\def\ed{\end{document}}
\def\de{\delta}

\def\r.BS{\rm reductive Borel-Serre compactification.\ }

\def\blem{\begin{lem}}
\def\bde{\begin{defn}}
\def\ede{\end{defn}}
\def\bex{\begin{ex}}
\def\eex{\end{ex}}
\def\bcl{\begin{claim}}
\def\ecl{\end{claim}}
\def\bcor{\begin{cor}}
\def\ecor{\end{cor}}

\def\bpr{\begin{prop}}
\def\epr{\end{prop}}
\def\elem{\end{lem}}
\def\bth{\begin{thm}}
\def\eth{\end{thm}}

\usepackage{tikz-cd}

\begin{document}
\title[Murmurations and Sato-Tate Law II: Beyond Riemann Hypothesis]{Murmurations\\ and Sato-Tate Conjectures\\ for High Rank Zetas of Elliptic Curves II:\\ Beyond Riemann Hypothesis}
\author{\bf Zhan Shi {\tiny and} Lin WENG}  
\date{}
\maketitle
\begin{abstract} As a continuation of our earlier paper \cite{SW},  we offer a new approach to murmurations and Sato-Tate laws for higher rank zetas of elliptic curves. 
Our approach here does not depend on the Riemann hypothesis for the so-called $a$-invariant $a_{E/\F_p;n}$ in rank $n(\geq 3)$ even for the Sato-Tate law, rather, on a much 
refined structure, a similar version of which was already observed earlier by Zagier and the senior author of this paper in \cite{WZ1}
when the rank $n$ Riemann hypothesis was established. Namely, instead of the rank $n$ Riemann hypothesis bounds 
$-1\leq \frac{a^{~^{~}}_{E/\mathbb F_p;n}}{2\sqrt{p^n}}\leq 1$ on which our first paper is based, we use  the asymptotic bounds
$-1\leq \frac{a^{~^{~}}_{E/\mathbb F_p;n}+(n-1)p+(n-5)}{2(n-1)\sqrt p}\leq 1$. Accordingly,  rank $n$ Sato-Tate law can be  established and 
rank $n$ murmurations can be formulated equally well, similar to the corresponding structures in the abelian framework for Artin zetas of elliptic curves.  
 \end{abstract}
\section{Introduction}
For a non-CM elliptic curve $\E$ over the field $\Q$ of rationals, distributions  of the associated $a$-invariants $a_{E/\F_p}$ defined by
$$a_{E/\F_p}:=1+p-\#E(\F_p),$$ or better, the associated arguments $\theta_{E/\F_p}$, defined by
 $$\cos\,\th_{E/\F_p}:=\frac{\ a^{~}_{E/\F_p}\ }{2\sqrt p},\qquad\th_{\E,p}\in[0,\pi]$$ for its $p$-reductions $E/\F_p$ satisfies the following famous Sato-Tate law: for any $0\leq \a<\b\leq\pi,$ 
$$ \lim_{N\to\infty}\frac{\#\Big\{p\leq N:p\ {\rm prime},\ \a\leq \th_{E/\F_p} \leq \b\Big\}}{\#\{p\leq N:p\ {\rm prime}\}}=\frac{\pi}{2}\int_\a^\b\sin^2\theta d\theta.$$

In \cite{SW}, for any fixed $n\in\Z_{\geq 1}$, we prove an analogue for the $a$-invariants $a_{E/\F_p;n}$ in rank $n$ of the $E/\F_p$'s
$$a_{E/\F_p;n}:=(p^n+1)-(p^n-1)\frac{\b_{E/\F_p;n}(0)}{\b_{E/\F_p;n-1}(0)}\qw 
\b_{E/\F_p;n}(0):=\sum_{\cV}\frac{1}{\#\Aut(\cV)}$$
where $\cV$ runs over all semi-stable vector bundles of rank $n$ and degree 0 on $E/\F_p$. 
To be more precisely, we have the following

\begin{thm} Let $\E/\Q$ be a non-CM elliptic curve. Denote its $p$-reduction by $E/\F_p$. Then, for $\alpha, \beta\in \R$ satisfying $0\le \alpha <\beta \le \pi$, we have
\bea
\lim_{N\rightarrow\infty}&\frac{\#\Big\{p\le N:p\ {\rm prime},\ \cos\alpha\ge\Delta_{E/\F_p,n}\ge\cos\beta\Big\}}{\#\{p\le N:p\ {\rm prime}\}}\\
=&\frac{2}{\pi}\int_\alpha^\beta\sin^2\theta d\theta.
\eea
Here the big $\Delta_{E/\F_{p},n}$ is defined by
\be
\Delta_{E/\F_{p},n}:=\bc
\sqrt{p}\cos{\theta_{E/\F_{p},2}}+\frac{1}{2}\Big(\sqrt{p}-\frac{1}{\sqrt{p}}\Big) & n=2
\\ 
&\\
\frac{\sqrt{p^{n-1}}}{n-1}\Big(\frac{\pi}{2}-\theta_{E/\F_{p},n}\Big)+\frac{1}{2}\sqrt{p}+\frac{1}{2}\frac{n-5}{(n-1)}\frac{1}{\sqrt{p}} & n\ge3 
\ec
\ee
and $\theta_{E/\F_{p},n}\in[0,\pi]$ is defined by
$$ \cos{\theta_{E/\F_{p},n}}:=\frac{\ a_{E/p;n}\ }{2\sqrt {p^n}},$$
\end{thm}
Our work in \cite{SW} heavily depends on the abelian Sato-Tate law, established by Taylor and his collaborators, based on the Riemann hypothesis for the zeta functions $\zeta_{E/\F_p;n}(s)$ in rank $n$. 

In this current work, we prove a new genuine Sato-Tate law for $a_{E/\F_p;n}$, which looks similar in appearance to the above rank $n$ Sato-Tate law, but based on a much strong estimate than the rank $n$ Riemann hypothesis used above.
\begin{thm}[Theorems 14, 15] Fix an integer $n\geq 3$. Let $\E/\Q$ be a non-CM elliptic curve. Denote its $p$-reduction by $E/\F_p$. Then we have
\bea
\lim_{N\to\infty}&\frac{\#\{p\leq N:p\ {\rm prime},\ \a\leq \Th_{E/\F_p;n}'\leq \b\}}{\#\{p\leq N: p: p\ {\rm prime}\}}\\
\lim_{N\to\infty}&\frac{\#\{p\leq N:p\ {\rm prime},\ \a\leq \Th_{E/\F_p;n}''\leq \b\}}{\#\{p\leq N: p: p\ {\rm prime}\}}\\
=&\frac{\pi}{2}\int_\a^\b\sin^2\theta d\theta.
\eea
Here, we set 
\bea
\cos\,\Th_{E/\F_p;n}':=&\De_{E/\F_p;n}'\quad{\rm{whenever}}\quad \left| \De_{E/\F_p;n}'\right|\leq 1\qquad and\\
\cos\,\Th_{E/\F_p;n}'':=& \De_{E/\F_p;n}''\quad{\rm{whenever}}\quad \left| \De_{E/\F_p;n}''\right|\leq 1\qquad{where}
\eea
\bea
\De_{E/\F_p;n}':=&\frac{1}{2(n-1)\sqrt p}\Big(a^{~}_{E/\F_p;n}-\big((5-n)-(n-1)p\Big),\qquad{and}\\
&\\
\De_{E/\F_p;n}'':=&\frac{a_{E/\F_{p};n}-\Big((5-n)+(n-1)a_{E/\F_p}-(n-1)p\Big)}{-6/\sqrt p}.\eea
\end{thm}

Thanks to the remarkable work \cite{HLOP}, we now have the murmuration structures in mathematics.  
Accordingly, for $n\geq 3$, in \cite{SW}, we introduce the murmuration functional 
\be
f_{r,n}(i):=\frac{1}{\#\cE_r[N_1,N_2]}\times\!\!\!\!\sum_{E\in\cE_r[N_1,N_2]}
\frac{1}{n-1}\cdot\big(a^{~}_{E/\F_{p_i},n}+(n-1)p_i+n-5\big).
\ee 
 as our non-abelian analogue of \cite{HLOP}. In this paper, we offer a new murmuration functional
$f_{r,n}^{\rm new}(i)$ by
setting
\bea
f_{r,n}^{\rm new}(i)=&\frac{1}{\#\cE_r[N_1,N_2]}\\
\times&\sum_{E\in\cE_r[N_1,N_2]}\Big(a^{~}_{E/\F_{p_i},n}+(n-1)p_i-(n-1)a_{E/\F_{p_i}}+(n-5)\Big)\frac{-p_i}{3}.\eea
Then with the same method as in \cite{SW}, we have the following
\begin{thm}[Theorem\,17]  For families of a regular (integral) elliptic curves $\E/\Q$'s, when plotting the points $(i,f_{r,n}^{\rm new}(i))$ $(i\ge 1, n\geq 3)$ in the sufficiently large range, the murmuration phenomenon appears in exactly the same way as that for the $(i,f_{r,1}(i))$'s.
\end{thm}

We end this introduction with the following illustrative figures, along with similar figures in \cite{SW}.
\begin{figure}[H]
    \centering
    \includegraphics[width=8.0cm]{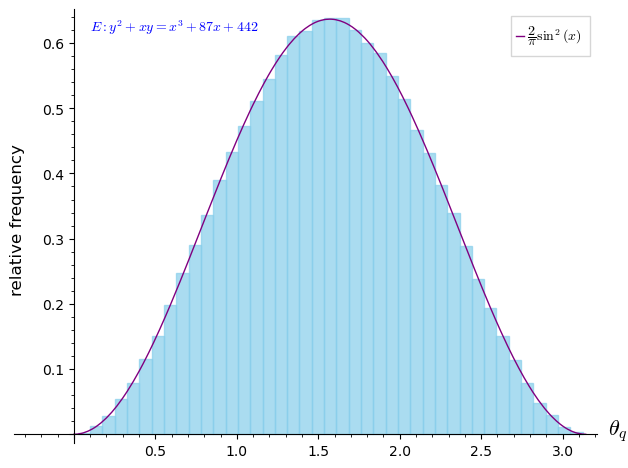}
    \caption{Sato-Tate distribution of rank 3 zeta function $\zeta_{E/\F_q,3}(s)$ in terms of $\De'$ over elliptic curve $\E/\Q: y^2 +xy= x^3 + 87 x +442$ and $q \le N =10,000,000$.}
    \label{hst}
\end{figure}
\begin{figure}[H]
    \centering
    \includegraphics[width=8.0cm]{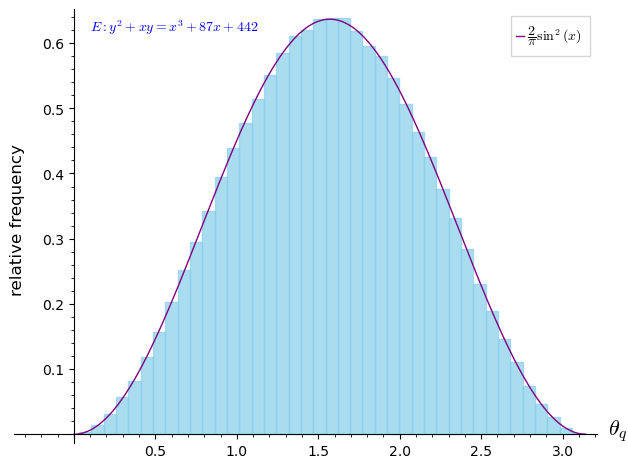}
    \caption{Sato-Tate distribution of rank 3 zeta function $\zeta_{E/\F_q,3}(s)$ in terms of $\De''$ over elliptic curve $\E/\Q: y^2 +xy= x^3 + 87 x +442$ and $q \le N =10,000,000$.}
    \label{hstnew}
\end{figure}\begin{figure}[H]
    \centering
    \includegraphics[width=8.0cm]{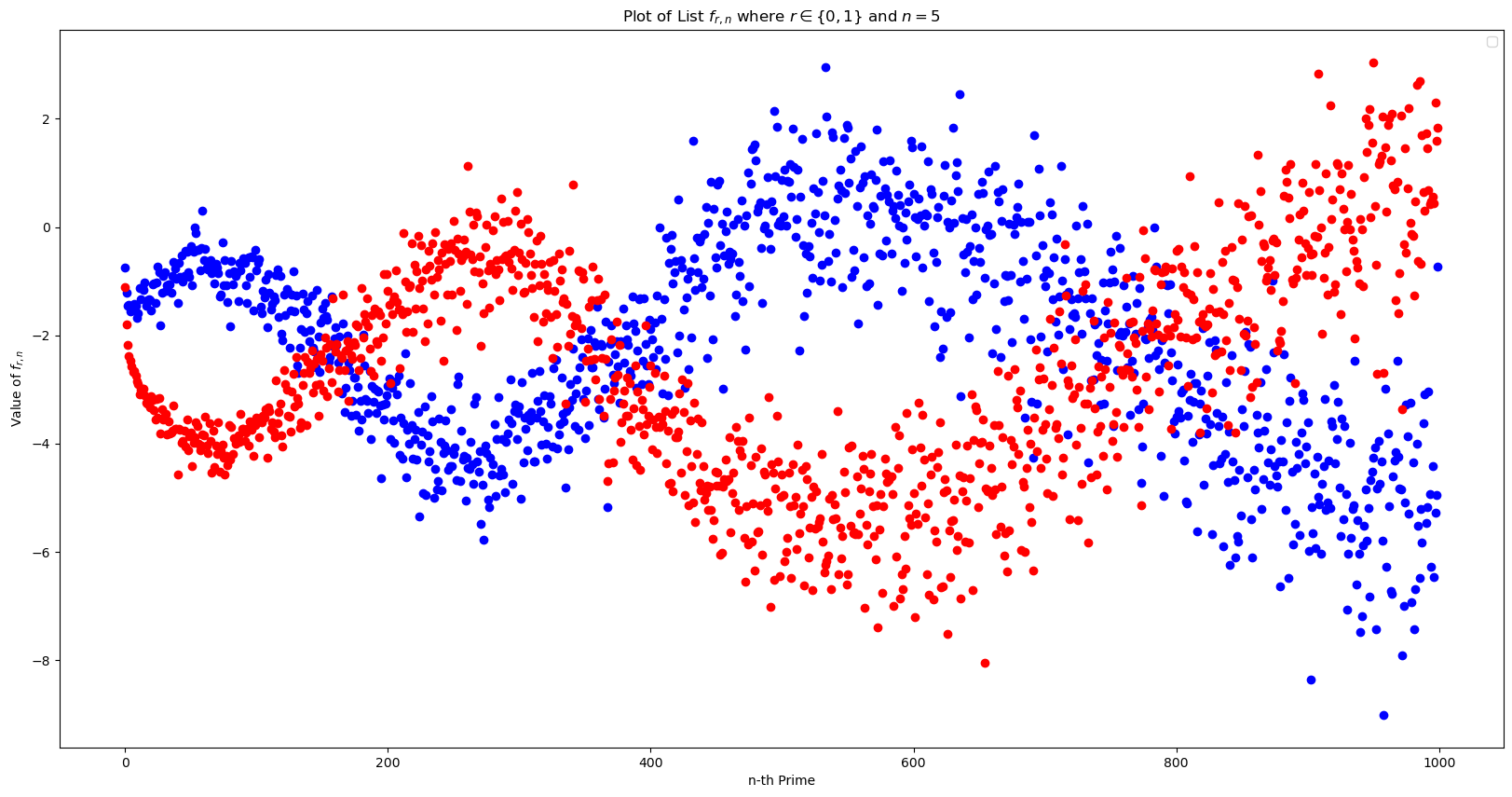}
    \caption{Plot of $f_{r,n}^{\rm new}(i)$ where $r\in{0,1}$ and $n=5$, for elliptic curves with conductor in $[7500,10000]$. $f_{0,n}^{\rm new}(i)$ is in blue and $f_{1,n}^{\rm new}(i)$ is in red.}
    \label{murnew0011}
\end{figure}
\begin{figure}[H]
    \centering
    \includegraphics[width=8.0cm]{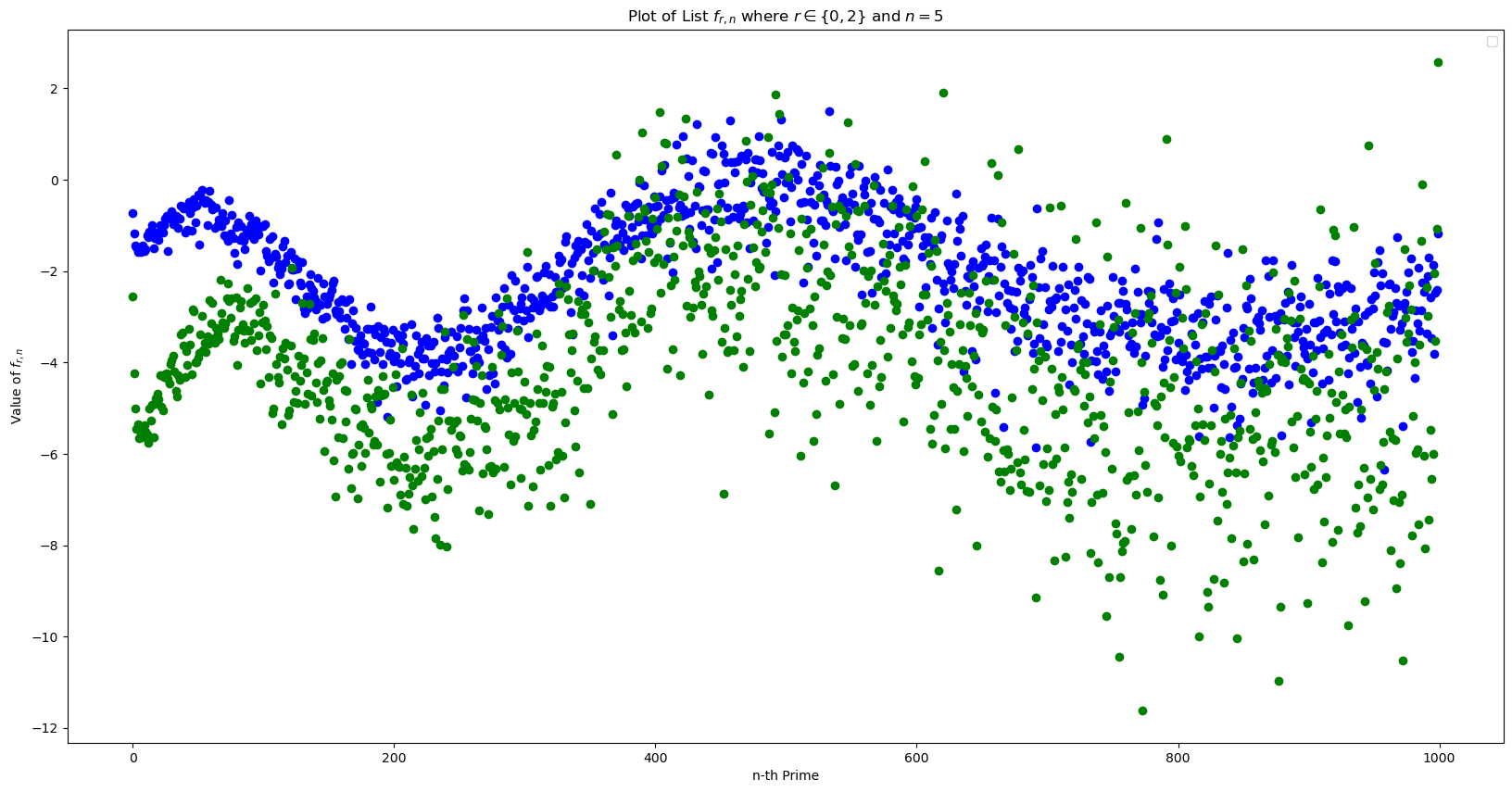}
    \caption{Plot of $f_{r,n}^{\rm new}(i)$ where $r\in{0,1}$ and $n=5$, for elliptic curves with conductor in $[7500,10000]$. $f_{0,n}^{\rm new}(i)$ is in blue and $f_{2,n}^{\rm new}(i)$ is in green.}
    \label{murnew0011}
\end{figure}

\section{Rank $n$ Sato-Tate Law based on Rank $n$ Riemann Hypothesis}
\subsection{Non-Abelian Zeta Function: Background}
By definition, for a(n integral regular projective) curve $X/\F_q$ (defined) over the finite field $\F_q$ with $q$ elements, its Artin zeta function is defined by
$$\zeta_{X/\F_q}:=\sum_{D\geq 0}\frac{1}{N(D)^s}\qquad\Re(s)>1,$$
where $D$ runs over all effective divisors on $X/\F_q$. Recall that, a formal sum $D=\sum_P n_PP$ with rational integral coefficients is called a divisor on $X/\F_q$, if, for almost all $P$, $n_P=0$. Moreover, such a $D$ is called effective, denoted by $D\geq 0$, if $n_P\in \Z_{\geq 0}$ for all $P$. As usual we set  $N(D):=\prod_PN(P)^{n_P}$. Here,  for each algebraic point $P$ on $X/\F_q$, if we set $k(P)$ to be the residue field of $X$ at $P$, and the norm $N(P)$ of $P$ is defined by $N(P):=q^{[k(P):\F_q]}$. 

From the standard zeta theory, it is well known that 

\begin{enumerate}
\item [(i)] (Rationality) $\zeta_{X/\F_q}$ is  rational function in $t=q^{-s}$. Indeed, there is a monic polynomial $P_{X/F_q}(t)\in \Z[t]$ of degree $2g$ such that
$$\wh \zeta_{X/\F_q}(s):=t^{g-1}\cdot \zeta_{X/\F_q}(s)=\frac{P_{X/F_q}(t)}{(1-t)(1-qt)t^{1-g}}$$
\item [(ii)] (Functional Equation)
$$\wh \zeta_{X/\F_q}(1-s)=\wh \zeta_{X/\F_q}(s).$$
\item [(iii)] (Geometric Interpretation of $\Res_{s=1}\wh \zeta_{X/\F_q}(s)$)
$$\Res_{s=1}\wh \zeta_{X/\F_q}(s)=\frac{\#\Pic(X)(\F_q)}{q-1}.$$
\end{enumerate}

Now let us regroup the $D$'s according to their rational equivalence classes $[D]$, then one arrives at
\bea
\zeta_{X/\F_q}(s)=&\sum_{[D]}\frac{\#\{D\geq 0:D\in [D]\}\smm\{0\}}{q-1}(q^{-s})^{\deg[D]}\\
=&\sum_{\cL\in\Pic(X/\F_q)}\frac{q^{h^0(X,\cL)}-1}{\#\Aut(\cL)}(q^{-s})^{\deg(\cL)}\eea where $\Pic(X/\F_q)$ denotes the Picard group of $X/\F_q$ and for each line bundle $\cL\in \Pic(X/\F_q)$, $h^0(X,\cL)$ denotes the dimension of the 0-th cohomology group $H^0(X,\cL)$ of $\cL$ over $X$. Indeed, this is because that there is an one-to-one correspondence between the set of  effective divisors $D$ in a rational equivalence class $[D]$ and  the set of divisors $(s)$ for a $\F_q$-line $\F_q\cdot s$ of a nontrivial global sections $s$ of the line bundle $\O_X(D)$.
 
As such, it is only temptation to define the rank $n$ zeta function for $X/\F_q$ by the formal summation
$\sum_{\cV}\frac{q^{h^0(X,\cV)}-1}{\#\Aut\cV}(q^{-s})^{\deg(\cV)}$
where $\cV$ runs over all rank $n$ vector bundles over $X/\F_q$.
However, this does not work due to the unboundedness of the the space of rank $n$ vector bundles with a fixed degree. One can sense this, in degree zero and rank 2, by examining the family
$\O_{\P^1}(-n)\oplus \O_{\P^1}(n)$ for $n\in\Z$ even over the projective line $\P^1/\F_q$. 

Fortunately, this unboundedness problem has already been tackled by geometers for decades. Following Mumford, to obtain a  nice family in fixed degree for rank $n$ vector bundles, stability condition should be included. Accordingly, one may try to introduce the rank $n$ zeta function for $X/\F_q$ by consider the  summation
\be
\sum_{\cV}\frac{q^{h^0(X,\cV)}-1}{\#\Aut\cV}(q^{-s})^{\deg(\cV)}\qquad\Re(s)>1,\ee
 where, instead over all rank $n$ vector bundles,  $\cV$ only runs over all rank $n$ semi-stable vector bundles over $X/\F_q$.
This works very nicely up to a certain reasonable stage. In fact in the senior author's first attempt \cite{W1}, this summation was used. In particular, by the standard zeta techniques, with helps from the vanishing theorem for semi-stable vector bundles, the Riemann-Roch and the duality, such defined series are verified to  enjoy some of the major standard zeta facts stated above, such as the rationality, the functional equation and the geometric interpretation of the residue at $s=1$.
 
However, examples suggest that, even with all semi-stable vector bundles counted, such defined generating functions do not satisfies the Riemann hypothesis. This indicates that there are something fundamentally wrong in this totality counting. It took the senior author several years to figure out the exact reason -- until a paper of Drinfeld (\cite{D}) on counting rank two super-cuspidal representations of the fundamental groups of $X/\F_q$ was introduced to him. 
\begin{defn}[\cite{W2}] Fix an integer $n\in\Z_{>0}$. Let $X/\F_q$ be an integral regular projective curve of genus $g$. Then  the rank $n$ zeta function of $X/\F_q$ is defined by
\bea
\wh \zeta_{X/\F_q;n}(s):=&(q^{-s})^{n(g-1)}\cdot  \zeta_{X/\F_q;n}(s)\\
:=&\sum_{\cV}\frac{q^{h^0(X,\cV)}-1}{\#\Aut\cV}(q^{-s})^{\chi(X,\cV)}\qquad(\Re(s)>1)\eea
where $\cV$ runs over all rank $n$ vector bundles over $X/\F_q$ \underline{whose degrees are} \underline{multiples of $n$}.
\end{defn}
Tautologically, as to be expected, by standard zeta techniques and the vanishing theorem for semi-stable vector bundles, the Riemann-Roch theorem and the Duality, we have
\begin{thm}[Zeta Facts \cite{W2}, see also \cite{WZ2}]\label{thm2} Fixed $n\in \Z_{\geq 1}$. The rank $n$ non-abelian zeta function 
$\zeta_{X/\F_q;n}(s)$ of an integral regular projective curve $X/\F_q$ satisfies the following standard zeta properties:
\begin{enumerate}
\item {\rm[(\underline{Naturality})]} We have
$$\zeta_{X/\F_q;1}(s)=\zeta_{X/\F_q}(s).$$
That is to say, the rank one zeta function $\zeta_{X/\F_q;1}(s)$ coincides with the classical Artin zeta function 
$\zeta_{X/\F_q}(s)$ of $X/\F_q$.
\item {\rm{ (\underline{Rationality})}} There exists a polynomial $P_{X/\F_q;n}(T)\in \Q[T]$ of degree $2g$, such that
$$
\zeta_{X/\F_q;n}(s)=\frac{P_{X/\F_q;n}(T)}{(1-T)(1-QT)}.
$$
In the above, we have set $T:=T_n:=t^n, Q:=Q_n:=q^n$.
\item {\rm (\underline{Functional Equation})} $\zeta_{X/\F_q;n}(s)$ satisfies the standard functional equation
$$
\wh \zeta_{X/\F_q;n}(1-s)=\wh \zeta_{X/\F_q;n}(s).
$$
\end{enumerate}
\end{thm}

In addition, we have the following
\begin{conj}[Riemann Hypothesis, \cite{W2}] The rank $n$ non-abelian zeta function 
$\zeta_{X/\F_q;n}(s)$ of an integral regular projective curve $X/\F_q$ satisfies the Riemann hypothesis. That is to say,
$$\zeta_{X/\F_q;n}(s)=0\Lra\Re(s)=\frac{1}{2}.$$
\end{conj}
This conjecture remains widely open, even its number theoretic analogue has been established (except when $n=1$ for  the lack of symmetry), up to a finite box depending on $n$ (\cite{W3}). Besides the classics, the first major breakthrough in this direction is the following:
  \begin{thm}[Weng-Zagier \cite{WZ1}]\label{thm3} Let $E/\F_q$ be an elliptic curve. Then, for $n\geq 2$,
 $\zeta_{E/\F_q,n}(s)$ satisfies the Riemann hypothesis.\footnote{We mention in passing that, besides this elliptic curve case, the Riemann hypothesis for $\zeta_{X/\F_q,n}(s)$ has been established  successfully  when
 \begin{enumerate} 
 \item[(0)] (Classical $n=1$) $X=E$ by Hasse, $X$ in general by Weil.
 \item[(i)] $n=2$ by H. Yoshida, see e.g. \S2 of arXiv:2201.03703.
 \item [(ii)] $n=3$ by Weng in \lq Riemann Hypothesis for Non-Abelian Zeta Functions of Curves over Finite Fields', arXiv:2201.03703.
 \item [(iii)] $g=2$ asymptotically by Shi, in preparation.
\end{enumerate}}
 \end{thm}

\subsection{Riemann Hypothesis in Rank $n$ for Elliptic Curves}

 By the rationality of the rank $n$ zeta functions for an elliptic curve $E/\F_q$,
 there exists a degree 2 polynomial $P_{E/\F_q;n}(T)\in \Q[T]$ such that
 $$\zeta_{E/\F_q;n}(s)=\frac{P_{E/\F_q;n}(T)}{(1-T)(1-QT)}.$$
 Moreover, if we define the so-called $\a$ and $\b$-invariants in rank $n$ for a curve $X/\F_q$ by
 $$\a_{E/\F_q;n}(d):=\sum_{\cV}\frac{q^{h^0(X,\cV)}-1}{\#\Aut(\cV)}\qqan \b_{E/\F_q;n}(d):=\sum_{\cV}\frac{1}{\#\Aut(\cV)}$$ where $\cV$ in the summations runs over all semi-stable vector bundles over $X/\F_q$ of rank $n$ and degree $d$. Clearly, $\b_{E/\F_q;n}(d)$ indeed counts semi-stable vector bundles
 naturally, by introducing the weight $\frac{1}{\#\Aut(\cV)}$ for each $\cV$,  being compatible with the language of algebraic stacks.
Much better, we have the following fundamental relation. 
 \begin{thm}[Counting Miracle. Theorem 3 of \cite{WZ1}] For all $n\geq 0$, we have
 $$\a_{E/\F_q;n+1}=\b_{E/\F_q;n}(0).$$
 \end{thm}
 Our approach is through a detailed analysis of the semi-stable vector bundles on $E/\F_q$ which then can be narrowly down to the so-called Atiyah bundles, together with a complicated combinatorial discussion.
As a direct  consequence, in (\cite{W2}), we show  the following
\begin{thm}[Equation 6 and Theorem 3 of \cite{WZ1}] With the same notation as above,
 $$P_{E/\F_q;n}(T)=\a_{E/\F_q;n}(0)\Big(1-a_{E/\F_q;n}T+Q_nT^2\Big)$$
 where the $a$-invariant of $E/\F_q$ in rank $n$ is defined by
 $$a_{E/\F_q;n}=(Q_n+1)-(Q_n-1)\frac{\b_{E/\F_q;n}(0)}{\b_{E/\F_q;n-1}(0)}.$$
 \end{thm}
 In particular, we have $\b_{E/\F_q;0}(0)=1$ and hence, when $n=1$,
 $$a_{E/\F_q;1}=(q+1)-(q-1)\frac{\b_{E/\F_q;1}(0)}{\b_{E/\F_q;0}(0)}=q+1-\#E(\F_q)=a_{E/\F_q}$$
 which is nothing but the classical $a$-invariant of $E/\F_q$. Consequently, the rank $n$-zeta function of $E/\F_q$ is completely determined by the $\b$-invariants.
 
 Accordingly, the Riemann hypothesis for the rank $n$ zeta function of $E/\F_q$ is equivalent to the fact that the degree two polynomial $1-a_{E/\F_q;n}T+Q_nT^2$ admits only non-real complex zeros. That is to say, the associated discriminant is strictly less than 0, or the same
 \be\label{eq6}\left|\frac{\ a_{E/\F_q;n}\ }{2\sqrt Q_n}\right|\leq 1.\ee
 With a sophisticated  combinatorial discussion, what we finally arrive in \cite{WZ1} is the following upper and lower bounds:
 \begin{thm}[Theorem 6 of \cite {WZ1}] For $n\geq 2$, we have
\be\label{eq7}1<\frac{\b_{E/\F_q;n}(0)}{\b_{E/\F_q;n-1}(0)}<\frac{\sqrt Q_n+1}{\sqrt Q_n-1}.\ee
\end{thm}
This then leads to the inequalities
\be\label{eq80}
2>a_{E/\F_q;n}>-2\sqrt Q_n
\ee
which are already noted in \cite{WZ1}. In other words, \eqref{eq7}, or the same \eqref{eq80}, is much refined than \eqref{eq6}. 

In fact, much refined structures on the $\b$-invariants in rank $n$ is exposed. 
\begin{thm}[Theorem 13 of \cite{WZ1}]\label{THM10} The $\b$-invariants in rank $n$ for elliptic curve $E/\F_q$ satisfies the following recursion formula: for $n\geq 1$,
$$(q^n-1)\b_{E/\F_q;n}(0)=(q^n+q^{n-1}-a_{E/\F_q})\b_{E/\F_q;n-1}(0)-(q^{n-1}-q)\b_{E/\F_q;n-2}(0),$$
with the initial conditions $\b_{E/\F_q;0}(0)=1$ and $\b_{E/\F_q;-1}(0)=0$.
\end{thm}
Consequently, all the $\b$-invariant invariants and hence the rank $n$ zeta function $\zeta_{E/\F_q;n}(s)$ are completely  determined  by $q$, $n$ and $a_{E/F_q}$. Based on this, we obtain the following asymptotic result
 
\begin{thm}[Theorem 6 of \cite{SW}] \label{thm6}
We have
\be
a_{E/\F_q,1}=a_{E/\F_q}, \qquad a_{E/\F_q,2}=1+a_{E/\F_q,1}-q,
\ee
and 
\be\label{eq10}
a_{E/\F_q,n}=(5-n)+(n-1)a_{E/\F_q,1}-(n-1)q+O\Big(\frac{1}{\sqrt{q}}\Big)\quad (n\geq 3)
\ee
In particular, for $n\geq 3$
\be\label{eq24}
a_{E/\F_q,n}\sim(5-n)+(n-1)a_{E/\F_q,1}-(n-1)q\ll 0\qquad (q\rightarrow\infty).
\ee
\end{thm}

Consequently, following the classical approach to formulate the Sato-Tate law for the distributions of the zeta zeros of elliptic curves $E/\F_p$'s associated to  $\E/\Q$, we are led to the construction of the big $\Delta$-distributions by single out $a_{E/\F_q}$ as stated in the Introduction. However, even it is very natural to use the Riemann hypothesis, or equivalently, the bounds
$-1\leq \frac{a_{E/\F_q;n}}{2\sqrt Q_n}\leq 1,$ to introduce $\theta_{E/\F_q;n}\in[0,\pi]$ via
$$\cos\, \theta_{E/\F_q;n}:=\frac{a_{E/\F_q;n}}{2\sqrt Q_n}$$  for an elliptic curve $\E/\Q$, one easily verify that the corresponding $\theta_{E/\F_{p_i};n}$'s have an obvious limit point $\frac{\pi}{2}$ when $n\geq 3$. 
This then yields the first structural distributions of the Dirac symbol $\de_{\pi/2}$ for the $\theta_{E/\F_{p_i};n}$'s. Unfortunately, $\theta_{E/\F_{p_i};n}-\frac{\pi}{2}$ is too tiny to be observed. Motivated by Theorem\,\ref{thm6}, a huge multiplicative factor $\sqrt{p_i^{n-1}}$ should be introduced so that the secondary level distributions of $\theta_{E/\F_{p_i};n}$ can be studied. However, with this enlargement, a further blow-up of additive scale $-(n-1)p_i$ is automatically introduced. It is for the purpose to eliminate this new complication, a term of $\frac{1}{2}\sqrt {p_i}$ is added, and hence to arrive finally at 
the normalized big $\Delta$-distributions:
\be
\Delta_{E/\F_{p},n}:=\bc
\sqrt{p}\cos{\theta_{E/\F_{p},2}}+\frac{1}{2}\Big(\sqrt{p}-\frac{1}{\sqrt{p}}\Big) & n=2
\\ 
&\\
\frac{\sqrt{p^{n-1}}}{n-1}\Big(\frac{\pi}{2}-\theta_{E/\F_{p},n}\Big)+\frac{1}{2}\sqrt{p}+\frac{1}{2}\frac{n-5}{(n-1)}\frac{1}{\sqrt{p}} & n\ge3. 
\ec
\ee
 In Theorem 4 of \cite{SW}, we are able to establish the following:
\begin{thm}[First Version of Sato-Tate Law in Rank $n$]\label{thm11} Fix a natural number $n\geq 2$. Let $\E/\Q$ be a non-CM  elliptic curve. For $\alpha, \beta\in \R$ satisfying $0\le \alpha <\beta \le \pi$, we have
$$
\lim_{N\rightarrow\infty}\frac{\#\{p\le N : p:\ {\rm prime},\ \cos\alpha\ge\Delta_{E/\F_p,n}\ge\cos\beta\}}{\#\{p\le N: p:\ {\rm prime}\}}=\frac{2}{\pi}\int_\alpha^\beta\sin^2\theta d\theta.
$$
\end{thm}
Our proof of this theorem is based on Taylor and his collaborators' works on the classical Sato-Tate law on the  abelian $a$-invariants $a_{E/\F_p}$'s.

\section{Sato-Tate Law in Rank $n$: Beyond Riemann Hypothesis}

\subsection{Riemann Hypothesis is Too Rough}

Recall that  \eqref{eq80} claims that 
$$
2>a_{E/\F_q;n}>-2\sqrt Q_n\qfo n\geq 2).
$$
This is clearly much strong than the bounds
$$\left|\frac{a_{E/\F_q;n}}{2\sqrt Q_n}\right|\leq 1,$$
which is well known to be equivalent to the Riemann hypothesis in rank $n$. 
This has already been observed in \cite{WZ1}.

Our first paper on the rank $n$ Sato-Tate law is based on the rank $n$ Riemann hypothesis. This leads to, as mentioned above,  a huge multiplicative factor $\sqrt{p^{n-1}}$, together with a further blowing up of additive scale $\sqrt p$, in the definition of the big normalized $\Delta$-distribution above.
Since this approach is guided by the Riemann hypothesis in rank $n$,  when such a big 
$\Delta$ was introduced in \cite{SW}, or much earlier in \cite{W3}, we thought, at the moment, that this approach was extremely natural.

However, recently, after examining the structures involved in the big $\De$ more carefully, we find out that
it is rather artificial to  use $$\frac{\sqrt{p^{n-1}}}{n-1}\Big(\frac{\pi}{2}-\theta_{E/\F_{p},n}\Big)$$
in the definition of $\De_{E/\F_p;n}$, since after all, $\Big(\frac{\pi}{2}-\theta_{E/\F_{p},n}\Big)$ is rather tiny, and hence should be of the same scale as $$\sin \Big(\frac{\pi}{2}-\theta_{E/\F_{p},n}\Big)=\cos\,\theta_{E/\F_{p},n}=\frac{a_{E/\F_p}}{2\sqrt {p^n}}.$$
Or equivalently
$$\frac{\sqrt{p^{n-1}}}{n-1}\Big(\frac{\pi}{2}-\theta_{E/\F_{p},n}\Big)=\frac{a_{E/\F_p}}{2\sqrt {p}(n-1)}.$$
which, in terms of Theorem\,\ref{thm11}, or better \eqref{eq10}, simply means that, up to normalization,  the big $\De$-distribution in rank $n\geq 3$ is essentially nothing but $\frac{\a_{E/\F_p}}{2\sqrt p}$.
Accordingly, for a genuine structure in  rank $n$, we are led to the following new yet more direct approach to the rank $n$ Sato-Tate law. 

For simplicity, in the sequel, unless otherwise is stated, let assume that $n\geq 3$.
Then, by \eqref{eq10}, we have as $q\to\infty$, 
$$a_{E/\F_q,n}=(5-n)+(n-1)a_{E/\F_q,1}-(n-1)q+O\Big(\frac{1}{\sqrt{q}}\Big).$$
Therefore, we set now
$$\De_{E/\F_p;n}':=\frac{1}{2(n-1)\sqrt p}\Big(a_{E/\F_q,n}+(n-1)p+(n-5)\Big)$$ and when
$\big|\De_{E/\F_p;n}'\big|\leq 1$,
we set $$\cos\,\Theta_{E/\F_p;n}'=\De_{E/\F_p;n}'.$$
The up-shot is the following:
\begin{thm}\label{thm13} Let $\E/\Q$ be a non-CM elliptic curve. Denote its $p$-reduction by $E/\F_p$. Then, for $\alpha, \beta\in \R$ satisfying $0\le \alpha <\beta \le \pi$, we have,
\bea
&\lim_{N\rightarrow\infty}\frac{\#\Big\{p\le N : p:\ {\rm prime},\ \alpha\le\Th_{E/\F_p,n}'\le\beta\Big\}}{\#\{p\le N: p:\ {\rm prime}\}}\\
=&\frac{2}{\pi}\int_\a^\b\sin^2\th d\theta.
\eea
\end{thm}
\bp
Note that, by Theorem\,\ref{thm6}, or better \eqref{eq10},  asymptotically, as $p\to\infty$,
$\big|\De_{E/\F_p;n}'\big|\leq 1$. Hence  by Theorem\,\ref{thm6}, or better \eqref{eq10} again, we have
\bea
&\lim_{N\rightarrow\infty}\frac{\#\Big\{p\le N : p:\ {\rm prime},\ \alpha\le\Th_{E/\F_p,n}'\le\beta\Big\}}{\#\{p\le N: p:\ {\rm prime}\}}\\
=&\lim_{N\rightarrow\infty}\frac{\#\Big\{p\le N : p:\ {\rm prime},\ \cos\alpha\ge\Delta_{E/\F_p,n}'\ge\cos\beta\Big\}}{\#\{p\le N: p:\ {\rm prime}\}}\\
=&\lim_{N\rightarrow\infty}\frac{\#\Big\{p\le N : p:\ {\rm prime},\ \cos\alpha\ge\frac{a_{E/\F_p}}{2\sqrt p}\ge\cos\beta\Big\}}{\#\{p\le N: p:\ {\rm prime}\}}\\
&(\text{by Theorem\,\ref{thm6}, or better \eqref{eq10}, based on the prime number theorem)}\\
=&\frac{2}{\pi}\int_\alpha^\beta\sin^2\theta d\theta\eea
Here in the last step, we have used the result of Taylor and his collaborators (\cite{BLGHT}, \cite{CHT}, \cite{HT}, \cite{T}) on the classical Sato-Tate law.
\ep
\subsection{Sato-Tate Law in Rank $n$: New Observing Spot}
After introducing the big $\De'$-distributions, our calculations indicates that in fact the difference
$$a_{E/\F_{p_i},n}-\Big((5-n)+(n-1)a_{E/\F_{p_i},1}-(n-1)p_i\Big)$$ for large $i$ is oscillating between
$\pm\frac{6}{\sqrt{p_i}}$. 

At the beginning, we thought that there might be a new type of distribution law
hidden behind this. However, through more detailed computations, we witness  that, numerically for all examples, 
$$\cos\,\Theta_{E/F_{q};n}'':=\frac{a_{E/\F_{q};n}-\Big((5-n)+(n-1)a_{E/\F_q}-(n-1)q\Big)}{-6/\sqrt q}$$
or better,  these new $\Theta_{E/F_{p_i};n}''$ obey the classical Sato-Tate law. 
In fact, we have the following
\begin{thm} Let $\E/\Q$ be a non-CM elliptic curve. Denote its $p$-reduction by $E/\F_p$. Then, for $\alpha, \beta\in \R$ satisfying $0\le \alpha <\beta \le \pi$, we have,

\bea
&\lim_{N\rightarrow\infty}\frac{\#\Big\{p\le N : p:\ {\rm prime},\ \alpha\le\Th_{E/\F_p,n}''\le\beta\Big\}}{\#\{p\le N: p:\ {\rm prime}\}}\\
=&\frac{2}{\pi}\int_\alpha^\beta\sin^2\theta d\theta\qfo n\geq 3).
\eea
\end{thm}
\bp We start with the following
\begin{thm}\label{thm15} Let $E/\F_q$ be an elliptic curve. Then for a fixed $n\geq 3$, we have
\be\label{eq180}
a_{E/\F_q;n}=(5-n)+(n-1)a_{E/\F_q}-(n-1)q-3\frac{a_{E/\F_q}}{q}+O\left(\frac{1}{q}\right).
\ee
\end{thm}
\bp
We use an induction on $n$.
Recall that from Theorem\,\ref{THM10},  
$$
(q^n-1)\b_{E/\F_q;n}(0)=(q^n+q^{n-1}-a_{E/\F_q})\b_{E/\F_q;n-1}(0)-(q^{n-1}-q)\b_{E/\F_q;n-2}(0),
$$
with the initial conditions $\b_{E/\F_q;0}(0)=1$ and $\b_{E/\F_q;-1}(0)=0$.
By definition, 
$$a_{E/\F_q;n}:=(q^n+1)-(q^n-1)\frac{\b_{E/\F_q;n}(0)}{\a_{E/\F_q;n}(0)}.$$
Thus, from the counting miracle that
$$\a_{E/\F_q;n}(0)=\b_{E/\F_q;n-1}(0),$$ we arrive at the structural recursion formula for $a_{E/\F_q;n}$:
\be\label{eq19}
a_{E/\F_q,n+1}
=1-q^n+a_{E/\F_q}+\frac{(q^n-q)(q^n-1)}{q^n+1-a_{E/\F_q,n}}
\ee
In particular,
\bea
a_{E/\F_q,2}=&1-q+a_{E/\F_q}\\
a_{E/\F_q,3}=&1-q^2+a_{E/\F_q}+\frac{(q^2-q)(q^2-1)}{q^2+1-a_{E/\F_q,2}}\\
=&1-q^2+a_{E/\F_q}+\frac{(q^2-q)(q^2-1)}{q^2+q-a_{E/\F_q}}\\
=&2+2a_{E/\F_q,1}-2q-3\frac{a_{E/\F_q}}{q}+O\left(\frac{1}{q}\right).\eea
This verifies \eqref{eq180} when $n=3$.

Now assume that
\eqref{eq180} holds for $n$. By \eqref{eq19}, we have

\bea
&a_{E/\F_q,n+1}\\
=&1-q^n+a_{E/\F_q}+\frac{(q^n-q)(q^n-1)}{q^n+1-a_{E/\F_q,n}}\\
=&1-q^n+a_{E/\F_q}+\frac{(q^n-q)(q^n-1)}{q^n+1-(5-n)-(n-1)a_{E/\F_q}+(n-1)q+3\frac{a_{E/\F_q}}{q}+O\left(\frac{1}{q}\right)}\\
=&(4-n)+na_{E/\F_q}-nq-3\frac{a_{E/\F_q}}{q}+O\left(\frac{1}{q}\right)
\eea
by a tedious, long yet trivial, calculation.
\ep
From Theorem\,\ref{thm15}, namely, \eqref{eq180},
we have
$$
a_{E/\F_q;n}=(5-n)+(n-1)a_{E/\F_q}-(n-1)q-3\frac{a_{E/\F_q}}{q}+O\left(\frac{1}{q}\right).
$$
Hence $$a_{E/\F_q;n}+(n-5)-(n-1)a_{E/\F_q}+(n-1)q=-3\frac{a_{E/\F_q}}{q}+O\left(\frac{1}{q}\right).$$
Therefore, asymptotically,
$$\frac{a_{E/\F_q;n}+(n-5)-(n-1)a_{E/\F_q}+(n-1)q}{-6/\sqrt q}=\frac{a_{E/\F_q}}{2\sqrt q}.$$
Therefore, our Theorem is equivalent to the classical Sato-Tate law for non-CM elliptic curve $\E/\Q$ established by Taylor and his collaborators (\cite{BLGHT}, \cite{CHT}, \cite{HT}, \cite{T}).
\ep

We end this section by the following comments: The Sato-Tate laws in higher ranks  can be observed from various spots, which have quite different geo-arithmetic meanings. These Sato-Tate laws exposes what should be called the secondary structures behind the Riemann hypothesis, and sometime are much refined  than the Riemann hypothesis in higher ranks. By contrast, all these Sato-Tate laws in higher ranks are essentially unique -- they are canonically equivalent to the canonical one.

\section{Rank $n$ Murmurations for Elliptic Curves}
To understand rank $n$ murmurations for elliptic curves $\E/\Q$'s of arithmetic rank $r$, in \cite{SW}, based on a weak version of \eqref{eq180} in Theorem\,\ref{thm15}, namely,
$$
a_{E/\F_q;n}=(5-n)+(n-1)a_{E/\F_q}-(n-1)q+O\left(\frac{a_{E/\F_q}}{q}\right),
$$
the following functional $f_{r,n}(i)$ on  the rank $n$ average value is introduced in \cite{SW}:
\bea
&f_{r,n}(i)\\
:=&\frac{1}{\#\cE_r[N_1,N_2]}\times\!\!\!\!\sum_{E\in\cE_r[N_1,N_2]}\bc a^{~}_{E/\F_{p_i},1}&n=1\\
a^{~}_{E/\F_{p_i},2}+p_i-1&n=2\\
\frac{1}{n-1}\cdot\big(a^{~}_{E/\F_{p_i},n}+(n-1)p_i+n-5\big)&n\geq 3\ec
\eea
where $N_1, N_2 \in \Z_{>0}$ satisfying $N_1\le N_2$, and $\cE_r[N_1,N_2]$ denotes the set of elliptic curve over $\Q$ of arithmetic rank $r$ with the conductor in the interval $[N_1, N_2]$.
Here, as in the rank one case,  for each isogeny class of elliptic curves $\E/\Q$, only a single representative elliptic curve is selected in $\cE_r[N_1,N_2]$. 

Next, we use Theorem\,\ref{thm15}, namely, a stronger estimate
$$
a_{E/\F_q;n}=(5-n)+(n-1)a_{E/\F_q}-(n-1)q-3\frac{a_{E/\F_q}}{q}+O\left(\frac{1}{q}\right),
$$
to introduce, for $n\geq 3$, a new murmuration functional
$f_{r,n}^{\rm new}(i)$ by
setting
\bea
f_{r,n}^{\rm new}(i)=&\frac{1}{\#\cE_r[N_1,N_2]}\\
\times&\sum_{E\in\cE_r[N_1,N_2]}\Big(a^{~}_{E/\F_{p_i},n}+(n-1)p_i-(n-1)a_{E/\F_{p_i}}+(n-5)\Big)\frac{-p_i}{3}.\eea
Then with the same method as in \cite{SW}, we have the following
\begin{thm}[Rank $n$ Murmurations]  For families of a regular (integral) elliptic curves $\E/\Q$'s, when plotting the points $(i,f_{r,n}^{\rm new}(i))$ $(i\ge 1, n\geq 3)$ in the sufficiently large range, the murmuration phenomenon appears in exactly the same way as that for the $(i,f_{r,1}(i))$'s.
\end{thm}
\bp
The same proof for Theorem \,4 of \cite{SW} works here as well since essentially, what $f_{r,n}^{\rm new}(i)$ really counts is $a_{E/\F_{p_i}}$ by Theorem\,\ref{thm15}.
\ep

\vskip 4.6cm
\noindent
Zhan SHI,\hskip 7.0cm Lin WENG,\\
shi.zhan.655@s.kyushu-u.ac.jp\hskip 3.60cm weng@math.kyushu-u.ac.jp\\
Graduate Program of Mathematics for Innovation\ \ Faculty of Mathematics\\ 
Kyushu University\hskip 5.50cm Kyushu University\\
Fukuoka, Japan\hskip 6.0cm Fukuoka, Japan
\end{document}